\documentclass[reqno,11pt]{amsart}
\usepackage{amsmath, latexsym, amsfonts, amssymb, amsthm, amscd}

\numberwithin{equation}{section}

\newtheorem{theorem}{Theorem}[section]
\newtheorem{lemma}[theorem]{Lemma}
\newtheorem{prop}[theorem]{Proposition}

\newtheorem{hypothesis}{Hypothesis}

\newcommand{\gga}{\gamma}            

\newcommand{\gep}{\varepsilon}       

\newcommand{\gl}{\lambda}

\newcommand{\go}{\omega}

\newcommand{\cB}{{\ensuremath{\mathcal B}} }

\newcommand{\cE}{{\ensuremath{\mathcal E}} }
\newcommand{\cH}{{\ensuremath{\mathcal H}} }

\newcommand{\cN}{{\ensuremath{\mathcal N}} }

\newcommand{\bbE}{{\ensuremath{\mathbb E}} }

\newcommand{\bbN}{{\ensuremath{\mathbb N}} }

\newcommand{\bbP}{{\ensuremath{\mathbb P}} }

\newcommand{\bbR}{{\ensuremath{\mathbb R}} }

\title{Convergence of approximations
of monotone gradient systems}

\author{Lorenzo Zambotti}
\address{Dipartimento di Matematica, Politecnico di Milano,
Piazza Leonardo da Vinci 32, 20133 Milano, Italy}
\email{lorenzo.zambotti\@@polimi.it}
\urladdr{http://www1.mate.polimi.it/\raisebox{0.11ex}{\tiny$\sim$}zambotti/}

\date{}

\begin{document}

\begin{abstract}
We consider stochastic differential equations in a Hilbert space, 
perturbed by the gradient of a convex potential. We investigate
the problem of convergence of a sequence of such processes.
We propose applications of this method to reflecting O.U. processes
in infinite dimension, to stochastic partial differential
equations with reflection of Cahn-Hilliard type and to interface models.
\end{abstract}

\maketitle

\section{Introduction}

Consider a separable Hilbert space $H$, which could be
finite or infinite dimensional, and a Stochastic Differential Inclusion of the form:
\begin{equation}\label{sde}
\left\{ \begin{array}{ll}
dX_t \in \left(AX_t - \partial U(X_t)\right)) \, dt +  dW_t,
\\ \\ X_0(x) \, = \, x\in H,
\end{array} \right.
\end{equation}
where $A$ is self-adjoint
in $H$, $U:H\mapsto ]-\infty,+\infty]$ is convex and lower semi-continuous 
with sub-differential $\partial U$ and $W$ a cylindrical white noise in $H$.

We recall that, if $U$ is Fr\'echet differentiable on $H$,
then $\partial U$ coincides with the gradient $\nabla U$. 
For general convex $U$, the subdifferential $\partial U(x)$ at $x\in H$ is
the set $\{y\in H: U(\xi)\geq U(x)+\langle y,\xi-x\rangle, \, \forall\,\xi\in H \}$
and the differential inclusion (\ref{sde}) is a formal way of writing which needs
to be made precise. Notice that $\partial U(x)$ is a non empty closed convex set
for all $x$ such that $U(x)<\infty$ and there exists an element of minimal norm $\partial_0 U(x)$.

If $\nabla U$ is Lipschitz-continuous, then existence and uniqueness of the SDE (\ref{sde})
is well known under suitable conditions on $A$. Moreover, in this case $X$
is reversible with respect to a probability measure $\nu:=\frac 1Z e^{-2U}d\mu$, where
$\mu$ is the Gaussian measure $\cN(0,(-2A)^{-1})$ and $X$ can be
characterized as the diffusion process in $H$
associated with the Dirichlet form obtained closing in $L^2(H,\nu)$
the symmetric bilinear form:
\begin{equation}\label{bifo}
\cE(\varphi,\psi) \, := \, \frac 12 \int\langle \nabla\varphi,\nabla\psi\rangle_H \, d\nu,
\qquad \varphi,\psi\in C^1_b(H).
\end{equation}

In this paper we consider a sequence of processes
$X^n$, solving (\ref{sde}) for some choice of $(H,A,U,W)=(H_n,A_n,U_n,W_n)$, and
we give conditions under which $X^n$ converges to a solution of (\ref{sde})
for some choice of $(H_\infty,A_\infty,U_\infty,W_\infty)$. More precisely, we
adopt here the Dirichlet form setting, and consider $X^n$ associated with $\cE^n$
for some choice of $(H_n,\nu^n)$ in (\ref{bifo}) and we give conditions for convergence of $X^n$ to
a process $X$ associated with $\cE$ for some choice of $(H,\nu)$. Our
result includes pointwise convergence of the transition semigroups $P^n$ of $X^n$
to the transition semigroups $P$ of $X$.

An important example of
this situation is the following: consider the Yosida approximations of $U$:
\begin{equation}\label{yo}
U_n(x) \, := \, \inf_{y\in H}\left\{U(y)+n\,\|x-y\|^2\right\}.
\end{equation}
Then it is well known that $U_n$ is smooth and $U_n$, respectively $\nabla U_n$,
converge to $U$, resp. $\partial_0 U$. An interesting particular case is:
\begin{equation}\label{example}
U(x) \, := \, \left\{ \begin{array}{ll}
0 \qquad & x\in K
\\ +\infty \qquad & x\notin K
\end{array} \right.
\end{equation}
where $K\subset H$ is a non-empty closed convex set. 
However we have several other interesting applications in mind besides Yosida
approximations: see the end of this introduction.

There is an extensive literature on equations of the type (\ref{sde}).
If $H=\bbR^d$ is finite-dimensional, then well-posedness and
existence of strong solutions of (\ref{sde})
even with more general drift and diffusion coefficients has been established by C\'epa in \cite{cepa}.
If dim $H=\infty$, then Da Prato and R\"ockner prove well-posedness 
of (\ref{sde}) in the class of weak solutions: see \cite{dpr}.
However, in the latter paper, $\nabla U$ is assumed to satisfy suitable integrability
conditions, and in particular the interesting case (\ref{example})
is not covered.

Notice also that (\ref{sde}) is naturally associated with the following second-order
elliptic operator:
\begin{equation}\label{Lu}
L\varphi(x) \, := \, \frac12 \, {\rm Tr}[D^2\varphi(x)] + \langle
Ax,\nabla\varphi(x)\rangle - \langle \partial_0 U(x),\nabla\varphi(x)\rangle,
\qquad \forall \, x\in D(A),
\end{equation}
where $\varphi$ is a smooth test-function. Da Prato, starting with the paper \cite{dp}, 
has investigated the analytical
properties of $L$ in a suitable $L^2(H,\mu)$ space, where the probability measure $\mu$
makes $L$ essentially self-adjoint. The same analytical approach is used by
Da Prato and Lunardi in \cite{dplu} in the finite-dimensional case. Again, if $H=\bbR^d$
then the paper \cite{dplu} covers essentially the most general situation, while, if
dim $H=\infty$, several interesting situations including (\ref{example})
are not covered by the existing literature.

All papers cited above (and many others) use the Yosida approximations
(\ref{yo}) and the approximating processes $X^n$, solving (\ref{sde})
with $U=U_n$. If $H=\bbR^d$, C\'epa \cite{cepa}
proves that $X^n$ converges almost surely to the solution of the limit differential inclusion;
moreover, Da Prato and Lunardi prove that the operator $L$ defined in (\ref{Lu}), with
a Neumann condition at the boundary of $K:=\{U<+\infty\}$, is
self-adjoint in $L^2(\bbR^d, e^{-2U}\, dx)$.
In infinite dimension, only weaker results are known under more restrictive assumptions.

Our approach allows to prove that under rather general conditions, the transition semigroup
of $X^n$
\[
P_t^n\varphi(x) \, := \, \bbE\left[\varphi(X^n_t(x))\right], \qquad \varphi\in C_b(H),
\ t\geq 0, \ x\in H,
\]
converges as $n\to\infty$ to a semigroup $P_t$, and in particular
that the finite dimensional distributions of the Markov process $X^n$ converge.
We can prove this, for instance, when $(-2A)^{-1}$ is trace-class,
for {\it any} convex lower semi-continuous 
$U:H\mapsto]-\infty,+\infty]$ such that $\mu(U<\infty)>0$ with $\mu:=\cN(0,(-2A)^{-1})$,
and $U_n$ is the Yosida approximation
of $U$. For instance, in the case (\ref{example}) with dim $H=\infty$ this result seems to be new.

We stress that we do not prove that the limit process $X$, which
exists by the Kolmogorov extension theorem, solves (\ref{sde}):
this can be (and has been) done in several interesting situations, like for
equations (\ref{cahi}) and (\ref{spder}) below, but only in particular
cases and using additional information about the model. For instance, a detailed description of the
right hand side of the integration by parts formula (\ref{ibpr}) below,
combined with the Fukushima decomposition, can lead to interesting
results: see \cite{za}.

\smallskip
We describe now briefly some applications of our general result.
First, we discuss the following Stochastic
Partial Differential Equation of Cahn-Hilliard type and reflection at $0$, that
has been considered in \cite{deza}:
\begin{equation}\label{cahi}
\frac{\partial u}{\partial t}=-
\frac{\partial^2 }{\partial \theta^2}
\left(\frac{\partial^2 u}{\partial \theta^2}
+\eta \right)  + \frac{\partial}{\partial\theta} \dot{W},
\qquad u\geq 0, \qquad d\eta\geq 0, \qquad \int u \, d\eta \, = \, 0
\end{equation}
where $u$ is a continuous function of $(t,\theta)\in [0,+\infty)\times[0,1]$,
$\eta$ is a locally finite positive measure on $(0,+\infty)\times[0,1]$ called the reflection
measure, preventing $u$ from becoming negative, and $\dot W$ is a space-time white noise 
on $[0,+\infty)\times[0,1]$. 
In this case we have $H=H^{-1}(0,1)$, Sobolev space of order $-1$, and $U$ of the form (\ref{example})
with $K:=\{x\in L^2(0,1): \, x\geq 0\}$. If we choose:
\[
U_n(x) \, := \, \left\{ \begin{array}{ll}
n\int_0^1 [x(\theta)\wedge 0]^2d\theta \qquad & x\in L^2(0,1)
\\ \\ +\infty \qquad & x\in H^{-1}(0,1)\backslash L^2(0,1)
\end{array} \right.
\]
then our method applies and we can prove convergence of $X^n$ as $n\to\infty$. Notice that $U_n$ is
not smooth in the topology of $H$ (in fact, $U_n$ is the Yosida approximation of $U$ in $L^2(0,1)$
rather than in $H$). We notice that $X^n$ is not a monotone sequence in $n$, and
no deterministic method is known to prove convergence of $X^n$.

Another interesting application is provided by interface models: in \cite{fuol}, Funaki and Olla
introduce a finite dimensional SDE whose solution $(\phi_i,i=1,\ldots,N)$ is stationary with 
respect to a 1-dimensional Gibbs measure. The process $\phi_i$ is always non-negative and
$\phi$ models the motion of an interface between  a gas and a liquid phase. Under a suitable
rescaling, the process $\phi$ converges in law to the solution a SPDE with reflection:
\[
\frac{\partial v}{\partial t} \, =  \, \frac{\partial^2
v}{\partial \theta^2} + \zeta + \dot{W},
\qquad 
v\geq 0, \qquad d\zeta\geq 0, \qquad \int v\, d\zeta=0.
\]
Also this problems fits in the general scheme we discuss here: a sequence $(\phi^N)$
associated with a Dirichlet form like (\ref{bifo}) with $\nu^n=\frac 1{Z^n} \, e^{-2U_n}\, dx$ and
$U_n$ convex, and a limit process $v$ associated with a Dirichlet form like (\ref{bifo})
with $\nu$ weak limit of $\nu^n$. Indeed, also the results of \cite{fuol}
can be proven using the techniques of this paper: we refer to \cite{zzz}.

\smallskip
The paper is organized as follows: in section 2 we present the setting and the
main results; in section 3 we show some applications; in section 4 we give a
tightness lemma and in section 5 we prove the main results.

\subsection{Notation}
For a real separable Hilbert space $H$,
we denote the scalar product by $\langle \cdot,\cdot\rangle_H$ and the associated norm 
by $\|\cdot\|_H$. If $J$ is Hilbert
space with scalar product  $\langle \cdot,\cdot\rangle_J$ and associated norm 
$\|\cdot\|_J$ such that $J$ is continuously embedded in $H$, 
we denote by $C_b(J)$, respectively $C_b^1(J)$, the space of all bounded
uniformly continuous functions on $J$, resp. bounded and uniformly
continuous together with the first Fr\'echet derivative. For all $\varphi\in C_b^1(J)$
we the directional derivative of $\varphi$ along $h\in J$ denote by $\partial_h  \varphi$:
\[
\partial_h  \varphi(x) \, := \, \lim_{t\to 0}\frac 1t(\varphi(x+th)-\varphi(x)),
\qquad x\in J.
\]
If $\nabla_J\varphi:J\mapsto J$ denotes the Fr\'echet
derivative of $\varphi$, then we have:
\[
\langle\nabla_J\varphi(x),h\rangle_J \, = \, \partial_h \varphi(x),
\qquad \forall \ x\in J, \ h\in J.
\]
If $\Gamma\subset H$
is closed and convex, we
denote by ${\rm Lip}(\Gamma)$ the set of all bounded $\varphi:\Gamma\mapsto\bbR$ with:
\[
[\varphi]_{{\rm Lip}(\Gamma)} \, := \, \sup \, \left\{
\frac{|\varphi(h)-\varphi(k)|}{\|h-k\|_H}, \
h\ne k, \, h,k\in\Gamma \right\} \, < \, \infty.
\]
Finally, if $D\subset H$ is dense,
we define ${\rm Exp}_D(H)\subset C_b(H)$ as the linear
span of $\{\cos(\langle h,\cdot\rangle_H)$, $\sin(\langle h,\cdot\rangle_H): h\in D\}$.

\section{Main result}

In this section we describe the general setting and the main
result. In the next section we show some concrete example.

All processes we consider take values in the fixed separable Hilbert
space $H$, which could be finite or infinite dimensional. In order to 
cover also convergence of finite dimensional approximations, we
consider a sequence $(H_n)_n$ of closed affine subspaces of $H$,
each itself a Hilbert space endowed with
a scalar product $\langle \cdot,\cdot\rangle_{H_n}$ and associated norm
$\|\cdot\|_{H_n}$, such that for a fixed constant $c>0$:
\begin{equation}\label{est}
\frac 1c \, \|h\|_H \,  \leq \, \|h\|_{H_n} \, \leq \,
c\, \|h\|_H, \qquad \forall \ h\in H_n, \ n\in\bbN.
\end{equation}
For all $n\in\bbN$ we consider a probability measure
$\nu^n$ on $H_n$ with topological support $K_n\subseteq H_n$ and we suppose
that $K_n$ is convex and closed. Moreover we consider
a continuous Markov process $X^n$ in $K_n$ such that:
\begin{hypothesis}\label{h} {\rm Let $c$ as in (\ref{est}). Then for
all $n\in\bbN$:
\begin{enumerate}
\item The transition semigroup $(P^n_t)_{t\geq 0}$ of $X^n$ 
acts on ${\rm Lip}(K_n)$ and for all $\varphi\in {\rm Lip}(K_n)$:
\begin{equation}\label{felle}
|P^n_t\varphi(x)-P^n_t\varphi(y)| \, \leq \, c \, 
[\varphi]_{{\rm Lip}(K_n)} \, \|x-y\|_H, \qquad x,y\in K_n,  \ t\geq 0.
\end{equation}
\item The following bilinear form is closable:
\[
{\cE}^{n}(\varphi,\psi) \, := \,
\frac 12 \, \int  \langle\nabla_{H_n} \varphi,\nabla_{H_n}\psi\rangle_{H_n} \, d\nu^n,
\qquad \forall \ \varphi,\psi\in C^1_b(H_n).
\]
and the closure $({\cE}^n,D({\cE}^n))$ is a Dirichlet form with associated
semigroup $(P^n_t)_{t\geq 0}$. 
\item For all $h$ in a dense subset $D_n\subset H_n$ there exists a finite signed
measure $\Sigma_h^n$ on $H_n$ such that:
\begin{equation}\label{sig0}
\int \partial_{h} \varphi\ d\nu^n \, = \, -
\int\varphi\ d\Sigma_h^n, \qquad \forall \ \varphi\in C_b^1(H_n).
\end{equation}
\end{enumerate}
}
\end{hypothesis}
\noindent
We recall that a {\it finite signed measure} on $H$ is a map $\Sigma$
from the Borel subsets $\cB(H)$ of $H$ to $\bbR$, such that $\Sigma(\emptyset)=0$ and
for any sequence $(E_h)_h\subset \cB(H)$ of pairwise disjoint sets:
\[
\Sigma\left(\bigcup_{h=0}^\infty E_h\right) \, = \, \sum_{h=0}^\infty \Sigma(E_h),
\]
where in the series in the right hand side we have absolute convergence. This notion
generalizes the definition of a measure to $\bbR$-valued set functions.
We also recall that we can associate to $\Sigma$ a finite positive measure $|\Sigma|$, called
the {\it total variation} of $\Sigma$, defined by:
\[
|\Sigma|(E) \, := \, \sup\left\{\sum_{h=0}^\infty |\Sigma(E_h)|: \,
E_h\in\cB(H) \ {\rm pairwise \ disjoint}, \ E=\bigcup_{h=0}^\infty E_h\right\}.
\]
Then $|\Sigma|$ is the smallest positive measure $\gga$ such that $|\Sigma(E)|\leq \gga(E)$
for all $E\in\cB(H)$. We refer to \cite[Chap. I]{afp}. 

Concerning the asymptotical behavior of $(H_n,\nu^n)_n$, we assume:
\begin{hypothesis}\label{hy} $ $ {\rm 
\begin{enumerate}
\item[(4)] $\nu^n$ converges weakly on $H$ to a probability
measure $\nu$, with convex topological support $K\subseteq H$
\item[(5)] For all $h$ in a dense subset $D\subset H$ 
there is a sequence $h_n\in D_n$ with $h_n\to h$ in $H$
such that for all $k_n\in H_n$ with $k_n\to k\in H$:
\begin{equation}\label{sig1.5}
\lim_{n\to\infty} \, \|h_n\|^2_{H_n} \, = \, \|h\|^2_H,
\qquad 
\lim_{n\to\infty} \, \langle h_n,k_n\rangle_{H_n}
\, = \, \langle h,k\rangle_H,
\end{equation}
and moreover there exist a finite signed
measure $\Sigma_h$ on $H$ and a sequence of
compact sets $(J_m)_m$ in $H$ such that:
\begin{equation}\label{sig1}
\lim_{n\to\infty}\int\varphi\, d\Sigma_{h_n}^n \, = \,
\int\varphi\, d\Sigma_h, \qquad \forall \ \varphi\in C_b(H),
\end{equation}
\begin{equation}\label{sig2}
\left|\Sigma_{h_n}^n\right|\left(H\backslash J_m\right) \, \leq \, \frac 1m.
\qquad \forall \ n,m\in\bbN.
\end{equation}
\end{enumerate}
}
\end{hypothesis}
\noindent
Notice that (\ref{sig2}) is a tightness condition for $\left|\Sigma_{h_n}^n\right|$.
For a sequence of probability measures, (\ref{sig1}) and (\ref{sig2}) are equivalent, but
this is not the case for signed measures: consider the example $H=\bbR$, $\Sigma^n=
\delta_{n+\frac 1n}-\delta_n$, with $\delta_a$ Dirac mass at $a$; then 
$\Sigma^n$ converges to the measure identically $0$ on $C_b(\bbR)$ (recall that
this is the space of all bounded {\it uniformly} continuous functions on $\bbR$),
but the sequence $|\Sigma^n|=\delta_{n+\frac 1n}+\delta_n$ is not tight. Therefore (\ref{sig2})
does not follow from (\ref{sig1}) and has to be proven separately.

By  (\ref{sig0}) and (\ref{sig1}), we have the following
integration by parts formula for $\nu$:
\begin{equation}\label{ibpr}
\int \partial_h \varphi \, d\nu \, = \, -\int \varphi
\, d\Sigma_h, \qquad \forall \ h\in D, \quad \forall
\ \varphi\in C^1_b(H),
\end{equation}
We denote the projection from $H$ to the element in $K_n$ with minimal distance
by $\Pi_{K_n}$:
\[
\Pi_{K_n}:H\mapsto K_n,  \qquad \|\Pi_{K_n}(x)-x\|_{H} \, \leq
\, \|k-x\|_{H} \qquad \forall \, k\in K_n.
\]
Under the above assumptions, we have the following Theorem, main result of the paper:
\begin{theorem}\label{main}
Suppose that Hypothesis \ref{h} and \ref{hy} hold. Then:
\begin{enumerate}
\item There exists a semigroup $(P_t)_{\geq 0}$ of operators
acting on $C_b(H)$ such that for all $\varphi\in C_b(H)$, $x\in K$,
setting $x_n:=\Pi_{K_n}(x)$:
\[
\lim_{n\to\infty} P_t^n\varphi\, (x_n) \, = \,
P_t\varphi(x), \qquad \forall \ t\geq 0.
\]
Moreover $[P_t\varphi]_{{\rm Lip}(H)}\leq c\,
[\varphi]_{{\rm Lip}(H)}$ for all $\varphi\in {\rm Lip}(H)$,
$t\geq 0$.
\item For all $x\in K$ there is a Markov process 
$(X_t)_{t\geq 0}$, defined on a probability space $(\Omega,\bbP_x)$, 
with state space $K$ and transition semigroup $(P_t)_{\geq 0}$, such that
$\bbP_x(X_0=x)=1$. Moreover $P_t\varphi(x)\to \varphi(x)$ as
$t\to 0$, for all $x\in K$, $\varphi\in C_b(H)$.
\item For all $\varphi_1,\ldots,\varphi_m\in C_b(H)$, $0\leq t_1\leq\ldots\leq t_m$
and $x\in K$, setting $x_n:=\Pi_{K_n}(x)$
\[
\lim_{n\to\infty}\,
\bbE_{x_n}\left[\varphi_1(X_{t_1}^n)\cdots\varphi_m(X_{t_m}^n)\right]\, = \, 
\bbE_x\left[\varphi_1(X_{t_1})\cdots\varphi_m(X_{t_m})\right].
\]
\item The following bilinear form is closable:
\[
{\cE}(\varphi,\psi) \, := \,
\frac 12 \, \int  \langle\nabla_H \varphi,\nabla_H\psi\rangle_H \, d\nu,
\qquad \forall \ \varphi,\psi\in {\rm Exp}_D(H).
\]
and the closure $({\cE},D({\cE}))$ is a Dirichlet form with associated
semigroup $(P_t)_{t\geq 0}$. Moreover ${\rm Lip}(H)\subset D(\cE)$
and $\cE(\varphi,\varphi)\leq[\varphi]^2_{{\rm Lip}(H)}$.
\item The stationary Markov process $\widehat X^n$ with transition
semigroup $P^n$ and initial distribution $\nu^n$
converges in law to the stationary Markov process $\widehat X$ with transition
semigroup $P$ and initial distribution $\nu$.
\end{enumerate}
\end{theorem}

\section{Applications}

Before proving Theorem \ref{main}, 
we discuss some interesting situations where the assumptions of Hypothesis \ref{h} and
\ref{hy} hold and therefore the results of Theorem \ref{main} apply.

\subsection{Gradient systems with convex potential}\label{gs}
Let now $H$ be of infinite dimension. We consider a closed
operator $A:D(A)\subset H\mapsto H$ such that:
\begin{enumerate}
\item $A$ is self-adjoint in $H$
and $\langle Ax,x\rangle_H \, \leq \, -\go\|x\|^2$ for all
$x\in D(A)$, with $\go>0$.
\item $Q:=(-2A)^{-1}$ is trace class in $H$ 
\end{enumerate}
Let $H_n=H$ for all $n\in\bbN$ and consider a convex lower semi-continuous function
$U:H\mapsto]-\infty,+\infty]$, setting $K:=\{U<+\infty\}$. We
define the Yosida approximations of $U$:
\[
U_n(x) \, := \, \inf_{y\in H}\left\{U(y)+n\,\|x-y\|_H^2\right\}.
\]
We recall that $U_n$ is convex, differentiable, and:
\begin{equation}\label{yosida}
U_n(x) \, \leq \, U_{n+1}(x),
\qquad \lim_{n\to\infty} U_n(x) \, = \, U(x), \qquad \forall \ x\in H.
\end{equation}
Moreover $\nabla U_n$ is Lipschitz-continuous with Lipschitz constant not larger than $2n$:
see \cite[Prop. 2.6, Prop. 2.11]{br}. By the Hahn-Banach Theorem, see
e.g. \cite[Proposition I.9]{bre},
and by the monotonicity of $n\mapsto U_n$, there exist constants $A\in\bbR$, $B>0$ such that:
\begin{equation}\label{bound}
U_n(x) \, \geq \, A \, - \, B\|x\|_H, \qquad \forall \ x\in
H, \quad n\in\bbN.
\end{equation}
Since $\nabla U_n$ is Lipschitz-continuous, there exists a unique strong solution $X^n(x)$ of:
\[
\left\{ \begin{array}{ll}
dX^{n}_t = (AX^{n}_t - \nabla_H U_{n}(X^{n}_t)) \, dt + dW_t,
\\ \\ X^n_0(x) \, = \, x\in H,
\end{array} \right.
\]
where $W$ is a cylindrical Wiener process in $H$.
Moreover:
\[
\|X^n_{t}(x)-X^n_{t}(y)\|  \, \leq \, 
\exp(-\go\, t) \, \|x-y\| , \qquad x,y\in H, \ t\geq 0,
\]
which yields (\ref{felle}).
We also define the probability measures on $H$: 
\[
\mu(dx) \, := \, \cN(0,Q)(dx), \qquad
\nu^n(dx) \, := \, \frac 1{Z^n} \,
\exp\left(-\, 2 \, U_{n}(x) \right)\, \mu(dx), 
\]
where $Z^n$ is a normalization constant. We also have for
all $h\in D(A)$:
\[
\Sigma_h^n (dx) \, := \, 2\, \left(\langle x,Ah\rangle_H -
\langle\nabla_H U_n(x),h\rangle_H\right)
\, \nu^n(dx).
\]
\begin{prop}\label{infini}
If $\mu(K)>0$,
then the assumptions of Hypothesis \ref{h} and \ref{hy} are fulfilled.
\end{prop}
\noindent
{\bf Proof}. Set $H_n:=H$, $D=D_n:=D(A)$. Since $\mu(K)>0$,
$\nu^n$ converges to:
\[
\nu \, := \, \frac 1Z \, e^{-2\, U} \, dx,
\]
where $Z>0$ is a normalization constant.
Let now $h=h_n\in D(A)$. For all $\varphi\in C^1_b(H)$:
\[
\exists \ \lim_n \int \varphi\, d\Sigma_h^n \, = \,
-\lim_n \int\partial_h\varphi\, d\nu^n \, = \, -\int\partial_h\varphi\, d\nu
\, =: \, \int \varphi\, d\Sigma_h.
\]
We do not know yet that $\Sigma_h$ is a signed measure. We set:
\[
|\Sigma_h|(H) \, := \,  \sup \, \left\{
\left|\int\varphi\, d\Sigma_h\right|: \, \varphi\in C^1_b(H),
\, \|\varphi\|_\infty\leq 1\right\}.
\]
Then we claim that $\Sigma_h^n$ and $\Sigma_h$ are indeed
finite signed measures on $H$ and:
\begin{equation}\label{cl}
\sup_n \, |\Sigma_h^n|(H) \, < \, \infty, \qquad
\lim_{n\to\infty} \ |\Sigma_h^n|(H) \, = \, |\Sigma_h|(H)
\, < \, \infty.
\end{equation}
Before proving (\ref{cl}) we show that it implies (\ref{sig1}) and (\ref{sig2}).
Indeed, setting $T_n:C_b(H)\mapsto\bbR$, $T_n\varphi:=\int
\varphi \, d\Sigma_h^n$, the family of linear functionals
$(T_n)_n$ is equicontinuous and converges on the dense
set $C^1_b(H)$ in the sup-norm and therefore converges on $C_b(H)$,
i.e. we have (\ref{sig1}). The proof of the density of $C^1_b(H)$ in $C_b(H)$ in the sup-norm
can be found in \cite[Theorem 2.2.1]{dpz3}.

On the other hand, consider a dense sequence $(x_j)_j$ in $H$ and set
$A_i^k:=\cup_{j=1}^n B(x_j,1/k)$. It is enough to prove that there exists $i$
such $|\Sigma_h^n|(H\backslash A_i^k)\leq 2^{-k}/m$: indeed, in this case 
$J_m:=\cap_k A_i^k$ is a compact set such that $|\Sigma_h^n|(H\backslash J_m)\leq 1/m$.
If we can not find such $i$, then for all $i$ there is $n(i)$ such that
$|\Sigma_h^{n(i)}|(A_i^k)\leq |\Sigma_h^{n(i)}|(H)-2^{-k}/m$. By the lower semi-continuity
of the total variation on open sets we find: 
$|\Sigma_h|(A_i^k)\leq |\Sigma_h|(H)-2^{-k}/m$
for all $i$ and therefore $|\Sigma_h|(H)\leq |\Sigma_h|(H)-2^{-k}/m$,
a contradiction. Therefore (\ref{sig2}) is proven.

\smallskip
We prove now (\ref{cl}). Let $Z$ be a $H$-valued random variable with 
distribution $\mu=\cN(0,Q)$ and $h\in D(A)$. 
We can suppose that $\langle -2Ah,h\rangle_H=1$. Setting:
$Y \, := \, Z \, + \, \langle Z,2Ah\rangle \, h$, we have that $Z=Y-\langle Z,2Ah\rangle h$
and $Y$ is independent of $\langle Z,2Ah\rangle$.
Notice that $Y\sim \cN(0,Q-h\otimes h)$ and $\langle Z,2Ah\rangle\sim\cN(0,1)$.
Then for all $z\in H$ we can write uniquely $z = y + t h$ with $t\in\bbR$ and 
$y\in (Ah)^\perp\subset H$ and with this notation:
\[
\mu(dz) \, = \, \cN(0,1)(dt) \cdot \cN(0,Q-h\otimes h)(dy).
\]
Let now $V:H\mapsto]-\infty,+\infty]$ be convex and lower semi-continuous.
Fix $y\in H$. Then by (\ref{bound}) the function $\bbR\ni t\mapsto v(t):=V(ht+y)+\frac{t^2}2$ is convex
and tends to $+\infty $ as $|t|\to\infty$. 
Necessarily $t\mapsto e^{-v(t)}$ is non-decreasing on a half-line
$(-\infty, t_0]$ and non-increasing on $(t_0,+\infty)$, where $v(t_0)=\min v$.
Then $e^{-v(\cdot)}$ has bounded variation on $\bbR$ and:
\[
\left|\frac{d}{dt} \, e^{-v(\cdot)}\right|(dt) \, = \, 
1_{(t\leq t_0)} \, \frac{d}{dt} \, e^{-v(\cdot)}(dt)
\, - \, 1_{(t> t_0)} \, \frac{d}{dt} \, e^{-v(\cdot)}(dt).
\]
We obtain that $\frac{d}{dt} \, e^{-v}$ is a finite signed measure on $\bbR$,
and:
\[
\int_\bbR \left|\frac{d}{dt} \, e^{-v}\right|(dt)
\, = \,  2\, e^{-v(t_0)} \, = \,  2\, e^{-\min \,v(\cdot)}.
\]
Now, by Fubini-Tonelli's Theorem, for all $\varphi\in C^1_b(H)$:
\begin{align*} &
\int\partial_h\varphi\, e^{-V} \, d\mu \, = \, \int \left[\int
\frac{\partial\varphi(ht+y)}{\partial t} \, e^{-V(ht+y)} \, \cN(0,1)(dt) \right]
\cN(0,Q-h\otimes h)(dy)
\\ & = \, \int \left[\frac1{\sqrt{2\pi}}
\, \frac{\partial}{\partial t} \, e^{-V(ht+y)-t^2/2}(dt) \right]\varphi(ht+y)\, 
\cN(0,Q-h\otimes h)(dy).
\end{align*}
Applying this formula to $V=2\, U_n$ and $V=2 \, U$ and
taking the supremum over all $\varphi\in C^1_b(H)$
such that $\|\varphi\|_\infty\leq 1$ we obtain:
\[
|\Sigma_h^n|(H) \, = \, \frac1{Z^n} \,
\int \sqrt{\frac2\pi}\, e^{- \min_{t\in\bbR} \, [2\, U_n(ht+y)+\frac{t^2}2]} \, \cN(0,Q-h\otimes h)(dy),
\]
\[
|\Sigma_h|(H) \, = \, \frac1{Z} \,
\int \sqrt{\frac2\pi}\, e^{- \min_{t\in\bbR} \, [2\, U(ht+y)+\frac{t^2}2]} \, \cN(0,Q-h\otimes h)(dy).
\]
Recall (\ref{bound}): it follows $2\, U_n(ht+y)+\frac{t^2}2\geq 2A - 2B^2\|h\|^2-2B\|y\|$
for all $t\in\bbR$, $y\in H$ and $n\in\bbN$ 
and therefore $|\Sigma_h^n|(H)+|\Sigma_h|(H)<+\infty$ for all $n\in\bbN$.

Now, fix $y\in H$ and set $v_n(t):=2\, U_n(ht+y)+\frac{t^2}2$, $v(t):=2\, U(ht+y)+\frac{t^2}2$.
The map $n\mapsto \min_t v_n(t)$ is monotone non-decreasing.
We want to prove that $\lim_n \min v_n=\min v$.
 
By (\ref{bound}): $\lim_{|t|\to\infty} v_n(t)=\infty$ for all $n\in\bbN$.
Let $c<\min v$: the set $I_n:=\{t:v_n(t)\leq c\}$ is
compact, $I_{n+1}\subset I_n$ and 
$\cap_n I_n=\emptyset$, since $t\in \cap_n I_n$ should
satisfy $v(t)=\lim_n v_n(t)\leq c<\min v$. By compactness,
there is $m\in\bbN$
such that $I_m=\emptyset$. Therefore for all $c<\min v$,
eventually $\min v_n>c$, i.e. $\lim_n \min v_n = \min v$.
By monotone convergence and since $Z^n\to Z$, (\ref{cl}) is proven. \quad $\square$

\medskip
In analogy with the finite dimensional example of the next subsection,
if $K\ne H$ then $\Sigma_h$ is expected to contain an infinite dimensional
boundary term: this has been explicitly computed in several
recent papers, for instance \cite{za}, \cite{boza}, \cite{fuis},
\cite{deza}. If $U$ is, as in (\ref{example}), equal to $0$ on
$K$ and to $+\infty$ on $H\backslash K$, then Theorem
\ref{main} gives a construction of the reflecting O.U.
process in $K$, in analogy with Fukushima's approach to
reflecting processes in finite dimension: see examples
4.4.2 and 4.5.3 in \cite{fot}.

For instance, let us consider the case $H=L^2(0,1)$, $A$ the
realization of $\partial^2$ with homogeneous Dirichlet
condition at $\{0,1\}$ and $U$ of the form (\ref{example})
with $K:=\{x\in L^2(0,1): x\geq -\alpha\}$
with $\alpha\geq 0$.
Then the Yosida approximation of $U$ is the distance from
$K$ in $H$, the process $X^n$ is the solution
of the SPDE:
\[
\frac{\partial X^n}{\partial t}=
\frac{\partial^2 X^n}{\partial \theta^2}
 + \frac{\partial^2 W}{\partial t\partial \theta} +
2n \left(X^n+\alpha\right)^-, \qquad \theta\in[0,1]
\]
where $W$ is a Brownian sheet, 
and the limit $X$ solves the following SPDE with reflection:
\begin{equation}\label{spder}
\frac{\partial X}{\partial t}=
\frac{\partial^2 X}{\partial \theta^2}
 + \frac{\partial^2 W}{\partial t\partial \theta} +
\eta(t,\theta), \qquad
X\geq -\alpha, \quad d\eta\geq 0, \quad
\int (X+\alpha)\, d\eta=0,
\end{equation}
where $\eta$ is a reflection term which prevents the continuous
solution $X$ from becoming less than $-\alpha$: see \cite{za}. 
If $\alpha>0$ then we are in the situation of Proposition \ref{infini}.
If $\alpha=0$ then $\mu(K)=0$ and Proposition \ref{infini} does
not apply directly. However in this case it is possible to
let first $n\to\infty$ with $\alpha>0$ and afterwards let
$\alpha\to 0$: both these limits satisfy the assumptions
of Hypothesis \ref{h} and \ref{hy}. 

Another interesting example is the Cahn-Hilliard equation with
reflection (\ref{cahi}) described in the introduction.

\subsection{Gradient systems with convex potential in finite dimension}\label{fd}

Let $H=H_n=\bbR^d$ and consider a convex lower semi-continuous
$U:\bbR^d\mapsto]-\infty,+\infty]$, such that 
$Z:=\int e^{-2U} \, dx<\infty$. By the convexity, the latter assumption is equivalent to:
\begin{equation}\label{infty}
\lim_{\|x\|\to\infty} U(x) \, = \, +\infty.
\end{equation}
We define $K:=\{U<\infty\}$. Then $K$ is a compact convex set 
and we assume that $K$ has non-empty interior. 
If we want to consider a case where $K\subset \bbR^d$ has empty interior, then 
we have to define $H'$ as the linear span of $K$ and work in $H'$ rather than
in $H$. Notice that $K$ necessarily has non-empty interior in $H'$.

As in the previous subsection, we define the Yosida approximations of $U$ as
$U_n(x) := \inf_{y\in {\bbR^d}}\left\{U(y)+n\,\|x-y\|^2\right\}$. Then
$U_n$ satisfies (\ref{yosida}) and (\ref{infty}). Let $W$ be 
a $d$-dimensional Brownian motion. Then by the
Lipschitz-continuity of $\nabla U_n$ there exists a unique strong solution $X^n(x)$ of: 
\begin{equation}\label{apfini}
dX^{n}_t = - \nabla U_{n}(X^{n}_t) \, dt + dW_t,
\qquad X^n_0(x) \, = \, x\in {\bbR^d},
\end{equation}
We also define the probability measure on $H$: 
\[
\nu^n(dx) \, := \, \frac 1{Z^n} \,
\exp\left(-\, 2 \, U_{n}(x) \right)\, dx, 
\]
where $Z^n$ is a normalization constant. Then, recalling (\ref{sig1}), we have
for all $h\in \bbR^d$:
\[
\Sigma_h^n (dx) \, := \, -2 \, \langle\nabla U_n(x),h\rangle
\, \nu^n(dx).
\]
\begin{prop}\label{fini}
In this setting, the assumptions of Hypothesis \ref{h} and \ref{hy} are fulfilled.
\end{prop}

\noindent
{\bf Proof}. The proof of Proposition \ref{infini} can
be repeated literally, substituting the Gaussian measure $\mu$ with the
Lebesgue measure $dx$ on $\bbR^d$ and the Gaussian measure $\cN(0,Q-h\otimes h)$
with the $(d-1)$-dimensional Lebesgue measure on $h^\perp$.
\quad $\square$

\medskip\noindent
In particular, Theorem \ref{main} applies to (\ref{apfini}). We notice that a stronger
result is proven in \cite{cepa}, namely almost sure convergence of $X^n$ as $n\to\infty$.
However the convergence of $\Sigma_h^n$ to $\Sigma_h$ in the sense
of (\ref{sig1}) and (\ref{sig2}), is quite general and of independent interest. 

Consider for instance the example (\ref{example}) with 
$K\subset\bbR^d$ a convex compact set with non-empty interior.
In this case $U_n=nd^2_K$, where $d_K(x)$ is the distance of $x$ from $K$.
Proposition \ref{fini} shows that (\ref{sig1}) and (\ref{sig2}) hold with:
\[
\Sigma^n_h(dx) \, = \, -4n\, d_K(x)\, \langle\nabla d_K(x),h\rangle \,
\frac 1{Z^n}\,  e^{-2n \, d_K^2(x)}\, dx,
\]
\[
\Sigma_h(dx) \, := \, - \frac1{|K|}\,
\langle \widehat n(x),h\rangle \, 1_{(x\in\partial K)}\cH^{d-1}(dx)
\]
as $n\to\infty$, where $\widehat n$ is an outer normal vector to the boundary of $K$, $\cH^{d-1}$
is the $(d-1)$-dimensional Hausdorff measure and $|K|\in(0,\infty)$ is the Lebesgue measure of $K$. 
When $K$ has
smooth boundary, then the normal vector and the boundary measure are classical objects; a
general bounded convex $K$ is a Lipschitz set, and 
the existence of the same objects follows from De Giorgi's work:
see e.g. \cite[Theorem 3.36]{afp}. The Dirichlet form approach to reflecting
Brownian motion in Lipschitz domains is developed by Fukushima in \cite{fot}
and Bass and Hsu in \cite{bahs}.

\subsection{SPDEs of Cahn-Hilliard type and Interface models}
Proposition \ref{infini} can be extended to other situations where
$U_n$ is not the Yosida approximation of $U$ but enjoys similar properties
For instance, this is the case for SPDEs with reflection of Cahn-Hilliard type,
considered in \cite{deza}, and interface models, as in \cite{fuol} and \cite{zzz}:
see the Introduction of this paper.

\section{Tightness of stationary approximations}

Let $\widehat X^n$ denote the stationary Markov process with transition
semigroup $P^n$ and initial distribution $\nu^n$. 
If $H=\bbR^d$, then in the following lemma
we can set $H^{(1)}=H^{(2)}=\bbR^d$. If $H$ is infinite dimensional,
then we fix Hilbert spaces $H^{(1)}$ and $H^{(2)}$ such that
$H$ is embedded in $H^{(1)}$ with Hilbert-Schmidt inclusion,
and $H^{(1)}$ is compactly embedded in $H^{(2)}$. We can suppose
that $H$ is dense in $H^{(2)}$.

\begin{lemma}\label{tightn}
For all $T>0$, the laws of $(\widehat X^{n})_{n\in\bbN}$ 
are tight in $C([0,T];H^{(2)})$.
\end{lemma}
\noindent{\bf Proof}. We claim that
for all $p>1$ there exists $C_p\in(0,\infty)$, independent of
$n$, such that:
\begin{equation}\label{estim}
\left( {\mathbb E}\left[
\left\|\widehat X^{n}_t-\widehat X^{n}_s\right\|^p_{H^{(1)}} \right]
\right)^{\frac1p} \, \leq \, C_p \, |t-s|^{\frac12}, \qquad
t,s\geq 0, \qquad  \forall \ n\in{\mathbb N}.
\end{equation}
To prove (\ref{estim}),
we fix $n\in\bbN$ and $T>0$ and use the Lyons-Zheng
decomposition, see e.g. \cite[Th. 5.7.1]{fot}, to write for
$t\in[0,T]$ and $h\in H_n$:
\[
\langle h,\widehat X^{n}_t - \widehat X^{n}_0\rangle_{H_n}
\, = \, \frac 12 \,
M_t \, - \, \frac 12 \, (N_T  - N_{T-t}), 
\]
where $M$, respectively $N$, is a martingale w.r.t. the natural
filtration of $\widehat X^{n}$, respectively of $(\widehat X^{n}_{T-t},
\ t\in[0,T])$.
Moreover, the quadratic variations are both equal to:
$\langle M\rangle_t =\langle N\rangle_t =t \cdot \|h\|^2_{H_n}$.
If $c_p$ is the optimal constant in the
Burkholder-Davis-Gundy inequality and $\kappa$ is the Hilbert-Schmidt norm of the
inclusion of $H$ into $H^{(1)}$, then (\ref{estim}) holds with $C_p=c_p \,
\kappa$.

Since the law of $\widehat X^{n}_0$ is
$\nu^n$ which converges as $n\to\infty$ in $H$ and
{\it a fortiori} in $H^{(2)}$,
tightness of the laws of $(\widehat X^{n})_{n\in\bbN}$
in $C([0,T];H^{(2)})$ follows
e.g. by Theorem 7.2 in Chap. 3 of \cite{ek}. \quad $\square$

\section{Convergence of the semigroups}

In this section we prove Theorem \ref{main}.
The proof is achieved using the theory of Dirichlet Forms, the uniform
Feller property (\ref{felle}) of $X^n$
and the integration by parts formula (\ref{ibpr}).
We define for all $\varphi\in C_b(K_n)$ the resolvent of $X^n$:
\[
R^n_\gl\varphi(x) \, := \, \int_0^\infty e^{-\lambda t} \,
P^n_t\varphi(x)\, dt, \qquad
x\in K_n, \ \lambda>0.
\]
We prove first
convergence of the resolvent $R^{n}$ in Proposition \ref{conres}
and then Theorem \ref{main}. 
\begin{prop}\label{conres} $ $
\begin{enumerate}
\item $(\cE,{\rm Exp}_D(H))$ is
closable in $L^2(\nu)$ and the closure
$({\cE},D({\cE}))$ is a symmetric Dirichlet form such that
${\rm Lip}(H)\subset D(\cE)$ and $\cE(\varphi,\varphi)\leq[\varphi]^2_{{\rm Lip}(H)}$.
The associated resolvent operator $(R_\gl)_{\gl>0}$ acts on $C_b(H)$. 
\item For all $\phi\in C_b(H)$ and
$x\in K$: $R^{n}_\gl\phi(\Pi_{K_n}(x))\to R_\gl\phi(x)$ as $n\to\infty$.
\end{enumerate}
\end{prop}
\noindent
We describe the idea of the proof: for $\varphi\in {\rm Lip}(H)$, by
the uniform Feller property (\ref{felle}) we have that $R^n_\gl\phi$
is an equicontinuous and equibounded family. By Arzel\`a-Ascoli's Theorem
we can extract converging subsequences on compact sets with large mass
with respect to $\nu^n$ and $\Sigma^n_h$. Now we consider the formula
which characterizes $R^n_\gl\phi$ for $\lambda>0$:
\[
\lambda \int   R^{n}_\gl\varphi \ \psi  \, d\nu^n  + \cE^n(R^{n}_\gl\varphi,
\, \psi)  \, = \, \int \psi \, \varphi \, d\nu^n, \qquad
\forall \ \psi\in D(\cE^n).
\]
We would like to pass to the limit, but $\cE^n$ contains the gradient of
$R^{n}_\gl\varphi$. However, if $\psi=\exp(i\langle \cdot,h\rangle)\in  {\rm Exp}_D(H)$, 
with $i^2=-1$, then we can use the integration
by parts formula (\ref{sig0}) and write:
\[
\cE^n(R^{n}_\gl\varphi,\, \psi) \, = \, -i\int R^n_\gl\varphi \ \psi  \, d\Sigma^n_h
\to -i\int F \, \psi  \, d\Sigma_h, \quad n\to\infty,
\]
where $F$ is a pointwise limit of $(R^n_\gl\varphi)_n$. Using (\ref{ibpr}),
the latter expression is equal to $\cE(F,\psi)$, i.e. we obtain:
\[
\lambda \int   F \, \psi  \, d\nu  + \cE(F,
\, \psi)  \, = \, \int \psi \, \varphi \, d\nu, \qquad
\forall \ \psi\in {\rm Exp}_D(H),
\]
which is very close to characterize $F$ as the $\gl$-resolvent of $\cE$
in $L^2(\nu)$ applied to $\varphi$. The proof of Proposition \ref{conres} makes these
arguments rigorous.

\smallskip
In the following proofs we use a number of times, often without further
mention, the following easily proven fact.
\begin{lemma}\label{easy} 
Let $E$ be a Polish space, $(M_n: \, n\in{\mathbb N}\cup\{\infty\})$ a sequence of
finite signed measures on $E$
and $(\varphi_n: \, n\in{\mathbb N}\cup\{\infty\})$ a sequence
of functions on $E$, such that: 
\begin{enumerate}
\item for all $\varphi\in C_b(E)$:
\[
\lim_{n\to\infty}\int\varphi\, dM_n \, = \,
\int\varphi\, dM_\infty
\]
\item there exists a sequence of
compact sets $(J_m)_m$ in $E$ such that:
\[
\lim_{m\to\infty} \ \sup_{n\in\bbN} \ 
\left|M_n\right|\left(E\backslash J_m\right) \, = \, 0.
\]
\item $(\varphi_n: \, n\in{\mathbb N}\cup\{\infty\})$ is 
equi-bounded and equi-continuous 
\item $\varphi_n$ converges pointwise to $\varphi_\infty$ on $\cup_m J_m$
\end{enumerate}
Then:
\[
\lim_{n\to\infty}\ \int_E \varphi_n \, dM_n \, = \, 
\int_E \varphi_\infty \, dM_\infty.
\]
\end{lemma}
\noindent
{\bf Proof}. We notice that by Arzel\`a-Ascoli's Theorem,
$\varphi_n$ converges uniformly to $\varphi$ on $J_m$
for all $m\in\bbN$. Moreover, by the Banach-Steinhaus Theorem
the norms of the functionals $C_b(E)\ni\varphi\mapsto \int_E \varphi\, dM_n$
are bounded, therefore $|M_n|(E)\leq C<\infty$ for all $n\in\bbN$. 
Then:
\begin{align*}
\left|\int_E \varphi_n \, dM_n \, -\, 
\int_E \varphi_\infty \, dM_\infty\right| \leq &
\left|\int_E
\left(\varphi_n - \varphi_\infty\right) \, dM_n\right|+
\left|\int_E \varphi_\infty \, \left(dM_n -dM_\infty\right)\right|
\end{align*}
and the second term in the right hand side tends to 0 by 
our first assumption. Now:
\begin{align*}
\left|\int_E \left(\varphi_n - \varphi_\infty\right) \, dM_n\right|
\leq & \int_{J_m} |\varphi_n - \varphi_\infty| \, d|M_n|+
\int_{E\backslash J_m} |\varphi_n - \varphi_\infty| \, d|M_n|
\\ \leq & \sup_{J_m}\, |\varphi_n - \varphi_\infty| \ C \, + \,
\|\varphi_n-\varphi_\infty\|_\infty \, |M_n|\left(E\backslash J_m\right).
\end{align*}
Taking the limsup as $n\to\infty$ and then letting $m\to\infty$ we
have the thesis. \quad $\square$

\medskip\noindent
{\bf Proof of Proposition \ref{conres}}. 
We divide the proof in several steps.

\medskip\noindent{\bf Step 1}. 
We recall that $\Pi_{K_n}:H\mapsto K_n$ is 1-Lipschitz in $H$, and therefore,
by (\ref{felle}):
\begin{equation}\label{jn}
\|(R^{n}_\gl\psi)\circ \Pi_{K_n}\|_\infty \, \leq \|\psi\|_\infty, \quad
[(R^{n}_\gl\psi)\circ \Pi_{K_n}]_{{\rm Lip}(H)} \, \leq  \, c \, [\psi]_{{\rm Lip}(H)},
\quad \forall \ \psi\in C_b^1(H).
\end{equation}
Fix $\psi\in C_b^1(H)$.
Let $(n_j)_j$ be any sequence in ${\mathbb N}$ and $(x_k)_k$ a 
countable dense set in $H$. With a diagonal procedure, we can find 
a subsequence $(m_i)_i$ and a function $F:\{x_k, k\in\bbN\}\mapsto {\mathbb R}$ 
such that $R^{n}_\gl\psi(\Pi_{K_n}(x_k))\to F(x_k)$ as $n=m_i\to\infty$
for all $k\in\bbN$. By (\ref{jn}), $F$ is Lipschitz on $\{x_k, k\in\bbN\}$
and therefore can be extended to a function in $\Psi_{\gl,\psi}\in{\rm Lip}(H)$
and:
\begin{equation}\label{subse}
\Psi_{\gl,\psi}(x) \, = \, \lim_{i\to\infty} R^{m_i}_\gl
\psi(\Pi_{K_{m_i}}(x)) \qquad \forall \, x\in H.
\end{equation}
Finally, by a diagonal procedure, we can suppose that such limit
holds along the same subsequence for all $\gl\in\bbN$. Notice that
in fact we are going to prove that the limit exists
as $n\to\infty$ for all $\gl>0$.
We define $\Delta:={\rm Span}\{\Psi_{\gl,\psi}: \psi\in C_b^1(H), \, \gl\in\bbN\}
\subset{\rm Lip}(H)$.

\medskip\noindent{\bf Step 2}. We would like to apply
the integration by parts formula (\ref{ibpr}) to
$\Psi_{\gl,\psi}$, which is however not in $C^1_b(H)$ but
only in ${\rm Lip}(H)$. However, notice that for all
$\varphi,\Phi\in C^1_b(H)$:
\begin{equation}\label{polariz}
\int \varphi \, \partial_h \Phi \, d\nu \, = \,
- \int \Phi \, \partial_h \varphi \, d\nu \, -\int \varphi
\, \Phi \, d\Sigma_h, \qquad \forall \ h\in D.
\end{equation}
If now $\Phi\in{\rm Lip}(H)$, then there exists a sequence
$(\Phi_m)_m\subset C_b^1(H)$ such that:
\[
\lim_m\, \Phi_m(x) \, = \, \Phi(x), \quad \forall \ x\in H,
\qquad \|\Phi_m\|_\infty+[\Phi_m]_{{\rm Lip}(H)}\, 
\leq \, \|\Phi\|_\infty+ [\Phi]_{{\rm Lip}(H)}.
\]
By (\ref{polariz}) we have that $\partial_h \Phi_m$
converges weakly in $L^2(\nu)$ to an element of $L^\infty(\nu)$ that
we call $\partial_h \Phi$ and with this definition
(\ref{ibpr}) holds for all $\Phi\in{\rm Lip}(H)$. Moreover,
we obtain in this way that $\nabla_H\Phi\in L^\infty(H,\nu;H)$
is well defined and:
\[
\cE(\Phi,\Phi) \, \leq \, \liminf_m \, \cE(\Phi_m,\Phi_m) \,
\leq \, \liminf_m \ [\Phi_m]_{{\rm Lip}(H)}^2 \, \leq \, [\Phi]_{{\rm
Lip}(H)}^2. 
\]
Moreover, for all $\Phi\in{\rm Lip}(H)$ it is possible to find a multi-sequence
$(\Phi_M)_M\subset {\rm Exp}_D(H)$, where $M=(m_1,\ldots,m_5)\in \bbN^5$, such that
$\Phi_M$ converges to $\Phi$ pointwise and:
\begin{equation}\label{den}
\sup_M\left(\|\Phi_M\|_\infty +
[\Phi_M]_{{\rm Lip}(H)} \right) \, < \, \infty,
\qquad \lim_{M}\, \cE(\Phi_M,\Psi) \, = \, \cE(\Phi,\Psi), 
\quad \forall \ \Psi\in{\rm Lip}(H),
\end{equation}
where $\lim_{M}$ means that we let
$m_1\to\infty$, then $m_2\to\infty$ and so on
until $m_5\to\infty$ (see \cite[Proposition 11.2.10]{dpz3} for
similar results).

\medskip\noindent
{\bf Step 3}. We want to prove now that for all $\gl\in\bbN$ and
$\Psi_{\gl,\psi}$ as in the first step:
\begin{equation}\label{st1}
\cE_\gl(\Psi_{\gl,\psi}, v)  \, := \,
\lambda \int   \,\Psi_{\gl,\psi}\, v \, 
d\nu + \cE(\Psi_{\gl,\psi}, v) 
\, = \, \int \psi \, v \, d\nu, \qquad
\forall \, v\in \Delta.
\end{equation}
First we prove (\ref{st1}) for $v\in {\rm Exp}_D(H)$.
Fix $h\in D$ and $h_n\in D_n$ as in Hypothesis \ref{hy} and set:
\[
\varphi_n(k) \, := \, \exp({i\langle h_n,\Pi_{H_n}k\rangle_{H_n} }), \qquad
\varphi(k) \, := \, \exp({i\langle h,k\rangle_H }), 
\qquad k\in H,
\]
where $i\in{\mathbb C}$ with $i^2=-1$ and $\Pi_{H_n}$ denotes 
the orthogonal projection from $H$ to $H_n$.
By Hypothesis \ref{hy}: $\|k-\Pi_{H_n}k\|_H\to 0$ for all $k\in H$.
Indeed, this is true for all $k\in D$ since there is a sequence $k_n\in H_n$
such that $k_n\to k$ and by density of $D$ in $H$ we conclude, since $\Pi_{H_n}$
is 1-Lipschitz continuous in $H$. Therefore, by (\ref{sig1.5}):
$\varphi_n(k)\to\varphi(k)$ for all $k\in H$.

Since $R^n_\gl$ is the resolvent operator associated with $\cE^n$:
\[
\cE_\gl^n(R^{n}_\gl\psi, \varphi_n)  \, := \,
\lambda \int   R^{n}_\gl\psi  \, \varphi_n \, d\nu^n  + \cE^n(R^{n}_\gl\psi, \varphi_n) 
\, = \, \int \psi \, \varphi_n \, d\nu^n.
\]
Notice that $\nabla_{H_n}\varphi_n=i\, h_n\, \varphi_n$. Then, by
the integration by parts formula \eqref{sig0}:
\[
2 \, \cE^n(R^{n}_\gl\psi, \varphi_n) \, = \, 
\int R^{n}_\gl\psi \ \|h_n\|^2_{H_n} \,
\varphi_n  \, d\nu^n \, - i \, \int R^{n}_\gl\psi \, \varphi_n\, d\Sigma_{h_n}^n.
\]
Since $\nu^n\rightharpoonup \nu$ and
$\Sigma_{h_n}^n\rightharpoonup \Sigma_h$
as $n\to\infty$, by (\ref{sig1.5}), (\ref{sig1}) and Lemma \ref{easy}: 
\[
\lim_{n\to\infty} \int g \ \|h_n\|^2_{H_n} \,
\varphi_n  \, d\nu^n \, - i \, \int g \, \varphi_n\, d\Sigma_{h_n}^n 
= \int g \, \|h\|^2_H \, \varphi \, d\nu -
i \, \int g \, \varphi\ d\Sigma_h, \
\forall g\in C_b(H).
\]
The crucial fact is now the following: by (\ref{sig1})-(\ref{sig2}),
(\ref{jn}) and Lemma \ref{easy}, we can substitute $g$ with
$R^n_\gl\psi$ in the last formula and prove that:
\begin{equation}\label{cru}
\lim_{i\to \infty} \int R^{m_i}_\gl\psi \, \varphi_{m_i} \,
d\Sigma^{m_i}_{h_{m_i}}\, = \, 
\int \Psi_{\gl,\psi}  \, \varphi \ d\Sigma_h.
\end{equation}
In particular we obtain:
\begin{align*}
\int \psi \, \varphi \, d\nu \, = \, & \lim_{i\to \infty}
\int \psi \, \varphi_{m_i} \, d\nu^{m_i} \, 
\, = \,  \int \Psi_{\gl,\psi}\left(\gl \, + \, \frac12 \,
\|h\|^2\right) \, \varphi \, d\nu
- i \, \frac 12 \int \Psi_{\gl,\psi} \, \varphi \ d\Sigma_h
\end{align*}
and by the integration by parts formula (\ref{ibpr}) the last expression is equal to
$\cE_\gl(\Psi_{\gl,\psi}, \varphi)$, i.e. we have proven (\ref{st1}) for
$v=\varphi$.
By linearity we obtain (\ref{st1}) for all $v\in{\rm Exp}_D(H)$. By
(\ref{den}) we obtain (\ref{st1}) for all $v\in \Delta$.

\medskip\noindent{\bf Step 4}. We want to prove now that
the bilinear form $(\cE,\Delta)$  is closable and the
closure is a Dirichlet form. 
By Lemma I.3.4 in \cite{maro}, it is enough 
to prove that if $(u_n)_n\subset \Delta$ and $u_n\to 0$ in
$L^2(\nu)$ then $\cE(u_n,v)\to 0$ for any $v\in \Delta$. 
But this is true, since, by (\ref{st1}), 
for all $u\in \Delta$ there exists some $\psi_u\in C_b(H)$ such that:
\[
\cE(u, v) \, = \, \int \psi_u \, v \, d\nu, \qquad
\forall \, v\in \Delta.
\]
We denote
the closure of $(\cE,\Delta)$ by $(\tilde\cE,D(\tilde\cE))$. We also obtain that 
$\Psi_{n,\psi}(x)=\tilde R_n\psi(x)$ for all $x\in K$ and
$n\in\bbN$, where $(\tilde R_{\gl})_{\gl>0}$ is the resolvent of
$\tilde\cE$.

\medskip\noindent{\bf Step 5}. Finally, we want to show that
$(\cE,{\rm Exp}_D(H))$ is closable and that the closure coincides with
$(\tilde\cE,D(\tilde\cE))$ constructed in the previous step.
To this aim we show first that $D(\tilde\cE)$
contains all Lipschitz functions on $K$ and in particular
${\rm Exp}_D(H)$. 

Consider $\psi\in {\rm Lip}(H)\subset D(\cE^{n})$:
by the general theory of Dirichlet Forms, 
\[
\psi\in D(\tilde{\cE}) 
\ \Longleftrightarrow \ 
\sup_{\lambda>0} \, \int
\lambda \, (\psi  -  \lambda \tilde R_\lambda\psi) \ \psi \
d\nu \ < \ \infty.
\]
By (\ref{felle}) we have:
\[
\int
\lambda \, (\psi -  \lambda R^{n}_\lambda\psi) \ \psi \
d\nu^{n}\ = \ {\cE}^{n}(\lambda R_\lambda^{n}\psi,\psi)
\ \leq \  [\psi]_{{\rm Lip}(H)}^2,
\]
so that letting $n\to\infty$:
\[
\int \lambda \, (\psi -  \lambda \tilde R_\lambda\psi) \ \psi \
d\nu \ \leq \ [\psi]_{{\rm Lip}(H)}^2,
\]
and therefore ${\rm Lip}(H) \subset D(\tilde{\cE})$. Since by construction
$\Delta\subset {\rm Lip}(H)$, then the closure of $(\cE,{\rm Lip}(H))$ is
$(\tilde\cE,D(\tilde\cE))$.
Now, in order to prove the density of ${\rm Exp}_D(H)$ in $D(\tilde{\cE})$,
we remark that the density with respect to the norm-topology
is equivalent to the density in the weak topology, which
follows from (\ref{den}) and from the density of ${\rm Lip}(H)$ in $D(\tilde{\cE})$.

Notice that the limit Dirichlet form $(\cE,D(\cE))$ does not depend
on the subsequences $(n_j)_j$ and $(m_i)_i$ chosen in step 1, since
it is the closure of $(\cE,{\rm Exp}_D(H))$.
Then $(\tilde\cE,D(\tilde\cE))=(\cE,D(\cE))$ is the Dirichlet form we wanted 
to construct and $R_\gl=\tilde R_\lambda$ is the associated resolvent operator.
In particular the limit in (\ref{subse}) does not
depend on the subsequence $(m_i)_i$ and
\[
R_{\gl}\psi(x) \, = \, \lim_{n\to\infty} \, R^{n}_\gl
\psi(\Pi_{K_n}(x)) \qquad \forall x\in K, \ \gl\in\bbN.
\]
We can now repeat the argument of step 1 and step 3 and obtain
that the latter formula holds for all $\gl>0$. 
$(\cE,D(\cE))$ is a Dirichlet Form, because
$R^{n}_\gl$ is given by a Markovian kernel, so that
$R_\gl$ is also Markovian and the result follows from Theorem 4.4 of
\cite{maro}. The Feller property follows from (\ref{felle}). 
By the density of ${\rm Lip}(H)$ in $C_b(H)$, $R^{n}_\gl
\psi(\Pi_{K_n}(x))$ converges to $R_{\gl}\psi(x)$ for all $\psi\in C_b(H)$.
\quad $\square$

\smallskip\noindent
{\bf Proof of Theorem \ref{main}}. We want to prove first that 
there exists a measurable kernel $(r_\gl(x,A): \, \gl>0,
x\in K, A\in{\cB}(H))$, such that for all $\varphi\in C_b(H)$ 
and $x\in K$:
\begin{equation}\label{ggg}
R_\gl \varphi(x) \, = \, \int \varphi(y) \, r_\gl(x,dy).
\end{equation}
By Lemma \ref{tightn} there exists a limit $\widehat X$
in distribution of $\widehat X^n$ along a subsequence $(n_i)_i$.
Then for all $\varphi,\psi\in C_b(H^{(2)})\subset C_b(H)$:
\begin{align*}
\int_0^\infty e^{-\gl \, t} \, \bbE\, [\psi(\widehat X_0) \, \varphi(\widehat X_t)]
\, dt \, = \,  & \lim_{i\to\infty} \int_0^\infty e^{-\gl \, t} \,
\bbE\, [\psi(\widehat X^{n_i}_0) \,
\varphi(\widehat X^{n_i}_t)] \, dt
\\ = \, & \lim_{i\to\infty} \int \psi \, R^{n_i}_\gl \varphi \
d\nu^{n_i} \, = \, \int \psi \, R_\gl \varphi \ d\nu.
\end{align*}
Since this is true for any $\psi\in C_b(H^{(2)})$ we obtain:
\[
R_\gl \varphi (x) \, = \, \int_0^\infty e^{-\gl \, t} \,
\bbE \, [\varphi(\widehat X_t) \, | \, \widehat X_0=x] \, dt , \qquad
\nu-{\rm a.e.} \ x.
\]
Therefore for $\gl>0$ fixed,
the measurable kernel $r_\gl(x,\cdot)$ exists for
$\nu$-a.e. $x$, in particular for $x$ in a set $K_0$ dense in $K$.
Let now any $x\in K\backslash K_0$ and $\varphi_n\in C_b(H)$ a
monotone sequence converging pointwise to $\varphi\in C_b(H)$
and let $x_m\in K$ converging to $x$. Then $R_\gl\varphi_n(x_m)$ 
converges to $R_\gl\varphi(x_m)$ by Dominated Convergence. Since $(R_\gl\varphi_n)_n$ is
equicontinuous, then necessarily also $R_\gl\varphi_n(x)$
must converge to $R_\gl\varphi(x)$, and the existence of
$r_\gl(x,\cdot)$ follows from general measure theory.

We prove now convergence of $P^n_t \varphi(\Pi_{K_n}(x))$ to $P_t \varphi(x)$. 
Let $\varphi\in{\rm Lip}(H)$:
by (\ref{felle}), for every $t>0$ and for any sequence $\gep_N\to 0$, 
we have pointwise convergence of $[P^n_t \varphi]\circ \Pi_{K_n}$ along 
a subsequence $(m_i)_i$ to $F\in C_b(H)$: see step 1 of the proof
of Proposition \ref{conres}. We want to
prove that $F=P_t \varphi$.
Let $m^{n,\varphi}$ be the spectral measure
of the generator of ${\cE}^{n}$ associated with $\varphi$: 
\[
\int_{-\infty}^0 \frac 1{(\lambda-x)^\ell} \, m^{n,\varphi}(dx) \, = \,
\int \varphi \left(R^{n}_\lambda\right)^\ell \varphi \ d\nu^n, \qquad \gl> 0,
\ \ell\in{\mathbb N}.
\]
Analogously, we denote the spectral measure of the generator of
${\cE}$ associated with $\varphi$ by $m^\varphi$. By (\ref{ggg}) and Lemma \ref{easy}, we obtain that
$\left(R^{n}_\lambda\right)^\ell$ converges pointwise to $\left(R_\lambda\right)^\ell$,
and therefore for all $\lambda>0$, $\ell\in{\mathbb N}$, by Proposition \ref{conres}:
\[
\int_{-\infty}^0 \frac 1{(\lambda-x)^\ell} \, m^{n,\varphi}(dx) \, 
\stackrel{n\to\infty}{\longrightarrow} \, 
\int_{-\infty}^0 \frac 1{(\lambda-x)^\ell} \, m^\varphi(dx).
\]
By the Stone-Weierstrass
theorem, the vector space spanned by the set of functions
$\{x\mapsto (\lambda-x)^{-\ell}, \, \lambda>0, \,
\ell\in{\mathbb N} \}$ 
is dense in the set $C_0((-\infty,0])$ of continuous
functions on $(-\infty,0]$ vanishing at $-\infty$. In
particular we obtain:
\[
\int \varphi \, P^{n}_t \varphi \ d\nu^n \, = \, \int_{-\infty}^0
e^{tx} \, m^{n,\varphi}(dx) \, \stackrel{n\to\infty} 
{\longrightarrow} \, \int_{-\infty}^0
e^{tx} \, m^\varphi(dx) \, = \, \int \varphi \, P_t \varphi \ d\nu.
\]
By polarization we have for all $\varphi,\psi\in C_b(H)$:
\[
\int \psi \, P^{n}_t \varphi \, d\nu^n \, \stackrel{n\to\infty}
{\longrightarrow} \, \int \psi \, P_t \varphi \, d\nu
\]
and we conclude that $\lim_{n\to\infty} P^{n}_t \varphi(\Pi_{K_n}(x)) =
P_t\varphi(x)$ for all $x\in K$. By the density of ${\rm Lip}(H)$ in $C_b(H)$, 
$P^{n}_t \varphi(\Pi_{K_n}(x))$ converges to $P_t\varphi(x)$ for all $\varphi\in C_b(H)$.

Arguing like for the existence of the kernel $r_\gl$, we can
prove that there exists a measurable kernel $(p_t(x,A): \, t\geq 0,
x\in K, A\in{\cB}(H))$, such that:
\[
P_t\varphi(x) \, = \, \int \varphi(y) \, p_t(x,dy),
\qquad \forall \ x\in K, \quad \varphi\in C_b(H).
\]
By Lemma \ref{easy}, we can now prove the convergence stated
in the third assertion of Theorem \ref{main}, obtaining also
that $(p_t(x,dy))_{t,x}$ satisfies the Chapman-Kolmogorov
equation. By the Kolmogorov extension Theorem we have existence
of the Markov process $X$. Convergence of $P_t\varphi(x)$ to
$\varphi(x)$ follows from the strong continuity of $P_t$ in
$L^2(\nu)$, the estimate $[P_t\varphi]_{{\rm Lip}(H)}\leq
c\, [\varphi]_{{\rm Lip}(H)}$ which follows from (\ref{felle}) and
the density of ${\rm Lip}(H)$ in $C_b(H)$.
The closability of $\cE$ has been
proven in Proposition \ref{conres}. The last assertion of Theorem \ref{main}
follows from Lemma \ref{tightn} and the third assertion.
\quad $\square$

\end{document}